\documentclass{article}    

\usepackage{amsmath,amssymb,graphicx,program}

\numberwithin{equation}{section}

\begin{document}           

\title{Formulae of numerical differentiation}

\author{Maxim Dvornikov
  \\
  Institute of Terrestrial Magnetism, Ionosphere and
  \\
  Radiowave Propagation (IZMIRAN)
  \\
  142190, Troitsk, Moscow region, Russia
  \\
  E-mail: maxdvo@izmiran.ru}         

\date{}

\maketitle                 

\thispagestyle{empty}

\begin{abstract}
We derived the formulae of central differentiation for the finding
of the first and second derivatives of functions given in discrete
points, with the number of points being arbitrary. The obtained
formulae for the derivative calculation do not require direct
construction of the interpolating polynomial. As an example of the
use of the developed method we calculated the first derivative of
the function having known analytical value of the derivative. The
result was examined in the limiting case of infinite number of
points. We studied the spectral characteristics of the weight
coefficients sequence of the numerical differentiation formulae.
The performed investigation enabled one to analyze the accuracy of
the numerical differentiation carried out with the use of the
developed technique.
\end{abstract}

\begin{center}
\end{center}


\section{Introduction}

In solving many mathematical and physical problems by means of
numerical methods one is often challenged to seek derivatives of
various functions given in discrete points. In such cases, when it
is difficult or impossible to take derivative of a function
analytically one resorts to numerical differentiation.

It should be noted that there exists a great deal of formulae and
techniques of numerical differentiation (see, for instance,
Ref.~\cite{DemMar63eng}). As a rule, the function in question
$f(x)$ is replaced with the easy-to-calculate function $\varphi
(x)$ and then it is approximately supposed that
$f^{\prime}(x)\approx\varphi^{\prime}(x)$. The derivatives of
higher orders are computed similarly. Therefore, in order to
obtain numerical value of the derivative of the considered
function it is necessary to indicate correctly the interpolating
function $\varphi (x)$. If the values of the function $f(x)$ are
known in $s$ discrete points, the function $\varphi (x)$ is
usually taken as the polynomial of $(s-1)$th power.

To find the derivative of functions having the intervals both
quick and slow variation quasi-uniform nets are used (see
Ref.~\cite{Kal78eng}). This method has an advantage since constant
small mesh width is unfavorable in this case, because it leads to
the strong enhancement of the function values table.

The problem of the numerical differentiation accuracy is also of
interest. The numerical differentiations formulae, taking into
account the values of the considered function both at $x>x_{0}$
and $x<x_{0}$ ($x_{0}$ is a point where the derivative is
computed), are called central differentiation formulae. For
instance, the formulae based on Stirling interpolating polynomial
can be included in this class. Such formulae are known to have
higher accuracy compared to the formulae, using unilateral values
of a function\footnote{The formulae using, for example, Newton
interpolating polynomial are attributed to this class of numerical
differentiation formulae.}, i.e., for instance, at $x>x_{0}$.

The range of numerical differentiation formulae based on different
interpolating polynomials is limited, as a rule, to finite points
of interpolation. All available formulae known at the present
moment are obtained for a certain concrete limited number of
interpolation points (see Refs.~\cite{AbrSti64,Bic41}). It can be
explained by the fact that the procedure of the finding of the
interpolating polynomial coefficients in the case of the arbitrary
number of interpolation points is quite awkward and requires
formidable calculations.

It is worth mentioning that the procedure of the numerical
differentiation is incorrect. Indeed, in Ref.~\cite{Kal78eng} it
was shown that it is possible to select such decreasing error of
the function in question which results in the unlimited growth of
the error in its first derivative.

Some recent publications devoted to the numerical differentiation
problem should be mentioned (see, e.g., Ref.~\cite{For88}). In
this work the finite difference formulae for real functions on one
dimensional grids with arbitrary spacing were considered.

The formulae of central differentiation for the finding of the
first and the second derivatives of the functions given in
$(2n+1)$ discrete points are derived in this paper. The number of
interpolation points is taken to be arbitrary. The obtained
formulae for the derivatives calculation do not require direct
construction of the interpolating polynomial. As an example of the
use of the developed method we calculate the first derivative of
the function $y(x)=\sin x$. The obtained result is studied in the
limiting case $n\to\infty$. We examine the spectral
characteristics of the weight coefficients sequence of the
numerical differentiation formulae for the different number of the
interpolation points. The performed analysis can be applied to the
studying of the accuracy of the numerical differentiation
technique developed in this work. It is found that the derived
formulae of numerical differentiation have a high accuracy in a
very wide range of spatial frequencies.

\section{Formulae for approximate values of
the first and the second derivatives} \label{f'f''deriv}

Without the restriction of generality we suppose that the
derivative is taken in the zero point, i.e. $x_{0}=0$. Let us
consider the function $f(x)$ given in equidistant points
$x_{m}=\pm mh$, where $m=0,\dots,n$ and $h$ is the constant value.
We can pass the interpolating polynomial of the $2n$th power
through these points
\begin{equation}
  \label{pol}
  P_{2n}(x)=\sum_{k=0}^{2n}c_{k}x^{k},
\end{equation}
the values of the function in points of interpolation
$f_{m}=f(x_{m})$ coinciding with the values of the interpolating
polynomial in these points: $P_{2n}(x_{m})=f_{m}$. Let us define
as $d_{m}$ the differences of the values of the function $f(x)$ in
diametrically opposite points $x_{m}$ and $x_{-m}$, i.e.
$d_{m}=f_{m}-f_{-m}$. We can present $d_{m}$ in the form
\begin{equation}
  \label{dm}
  d_{m}=2\sum_{k=0}^{n-1}c_{2k+1}h^{2k+1}m^{2k+1}.
\end{equation}
To find the coefficients $c_{2k+1}$, $k=0,\dots n-1$, we have
gotten the system of inhomogeneous linear equations with the given
free terms $d_{m}$. It will be shown below that this system has
the single solution.

We will seek the solution of the system [Eq.~\eqref{dm}] in the
following way
\begin{equation}
  \label{c2k}
  c_{2k+1}=\frac{1}{2h^{2k+1}}
  \sum_{m=1}^{n}d_{m}\alpha_{m}^{(2k+1)}(n),
\end{equation}
where $\alpha_{m}^{(2k+1)}(n)$ are the undetermined coefficients
satisfying the condition
\begin{equation}
  \label{alpha}
  \sum_{m=1}^{n}\alpha_{m}^{(2l+1)}(n)m^{2k+1}=\delta_{lk},
  \quad l,k=0,\dots,n-1.
\end{equation}
Thus, the system of equations [Eq.~\eqref{dm}] is reduced to the
equivalent, but more simple system [Eq.~\eqref{alpha}], in which
for each fixed number $k=0,\dots,n-1$ it is necessary to find the
coefficients $\alpha_{m}^{(2l+1)}(n)$.

Let us resolve the system of equations [Eq.~\eqref{alpha}]
according to the Cramer's rule:
\begin{equation}
  \label{kram}
  \alpha_{m}^{(2l+1)}(n)=
  \frac{\Delta_{m}^{(2l+1)}(n)}{\Delta_{0}(n)},
\end{equation}
where
\begin{equation}
  \label{delta0}
  \Delta_{0}(n)=
  \begin{vmatrix}
  1 & 2 & \dots & n
  \\
  1 & 2^{3} & \dots & n^{3}
  \\
  \hdotsfor{4}
  \\
  1 & 2^{2n-1} & \dots & n^{2n-1}
  \\
  \end{vmatrix}=
  n!\prod_{1\leq i<j\leq n}(j^{2}-i^{2}){\not=}0,
\end{equation}
\[
  \Delta_{m}^{(2l+1)}(n)=
  \begin{vmatrix}
  1 & 2 & \dots & m-1 & 0 & m+1 & \dots & n
  \\
  \hdotsfor{8}
  \\
  1 & 2^{2l+1} & \dots & (m-1)^{2l+1} & 1
  & (m+1)^{2l+1} & \dots & n^{2l+1}
  \\
  \hdotsfor{8}
  \\
  1 & 2^{2n-1} & \dots & (m-1)^{2n-1} & 0
  & (m+1)^{2n-1} & \dots & n^{2n-1}
  \end{vmatrix}.
\]
In Eq.~\eqref{delta0} we used the formula for the calculation of
the Vandermonde determinant. From Eq.~\eqref{delta0} it follows
that the determinant of the system of equations
[Eq.~\eqref{alpha}] is not equal to zero, i.e. the system of
equations [Eq.~\eqref{dm}] has the single solution.

The most simple expression for $\Delta_{m}^{(2l+1)}(n)$ is
obtained in the case of $l=0$ that corresponds to a calculation of
the first-order derivative
\begin{equation}
  \label{deltam1}
  \Delta_{m}^{(1)}(n)=(-1)^{m+1}
  \left(
  \frac{n!}{m}
  \right)^{3}
  \prod_{
  \substack{1\leq i<j\leq n
  \\
  i,j{\not=}m}}(j^{2}-i^{2}).
\end{equation}
From Eq.~\eqref{kram} as well as taking into account
Eqs.~\eqref{delta0} and \eqref{deltam1} we get the expression for
the coefficients $\alpha_{m}^{(1)}(n)$
\begin{equation}
  \label{alpha1}
  \alpha_{m}^{(1)}(n)=\frac{1}{m\pi_{m}(n)},
\end{equation}
where
\begin{equation}
  \label{pim}
  \pi_{m}(n)=\prod_{
  \substack{k=1
  \\
  k{\not=}m}}^{n}
  \left(
  1-\frac{m^{2}}{k^{2}}
  \right).
\end{equation}
It should be noted that one can similarly get the expression for
the coefficients $\alpha_{m}^{(2n-1)}(n)$ which is presented in
the following way
\[
  \alpha_{m}^{(2n-1)}(n)=\frac{(-1)^{n+1}m}{(n!)^{2}\pi_{m}(n)}.
\]

Taking into account Eqs.~\eqref{pol}-\eqref{c2k} we finally get
the formula for the first derivative of the function $f(x)$
\begin{equation}
  \label{f'}
  f^{\prime}(0)\approx P_{2n}^{\prime}(0)=
  \frac{1}{2h}
  \sum_{m=1}^{n}\alpha_{m}^{(1)}(n)(f_{m}-f_{-m}).
\end{equation}
The algorithm for the computation of the coefficients
$\alpha_{m}^{(1)}(n)$ is presented in the appendix~\ref{append}
and the results for the certain concrete number of the
interpolation points $(2n+1)$ are given in Tab.~\ref{alphatab}.
\begin{table}
\begin{center}
  \begin{tabular}{|l||c|c|c|c|c|c|} \hline
    $n$ & $\alpha_{1}^{(1)}(n)$ & $\alpha_{2}^{(1)}(n)$ &
    $\alpha_{3}^{(1)}(n)$ & $\alpha_{4}^{(1)}(n)$ &
    $\alpha_{5}^{(1)}(n)$ & $\alpha_{6}^{(1)}(n)$ \\ \hline \hline
    $1$ & $1$ & $0$ & $0$ & $0$ & $0$ & $0$ \\ \hline
    $2$ & ${4}/{3}$ & $-{1}/{6}$ & $0$ & $0$ & $0$ & $0$ \\ \hline
    $3$ & ${3}/{2}$ & $-{3}/{10}$ & ${1}/{30}$ & $0$ & $0$ & $0$ \\ \hline
    $4$ & ${8}/{5}$ & $-{2}/{5}$ &
    ${8}/{105}$ & $-{1}/{140}$ & $0$ & $0$  \\ \hline
    $5$ & ${5}/{3}$ & $-{10}/{21}$ &
    ${5}/{42}$ & $-{5}/{252}$ & ${1}/{630}$ & $0$  \\ \hline
    $6$ & ${12}/{7}$ & $-{15}/{28}$ &
    ${10}/{63}$ & $-{1}/{28}$ & ${2}/{385}$ & $-{1}/{2772}$  \\ \hline
  \end{tabular}
  \caption{\label{alphatab} Values of the coefficients
  $\alpha_{m}^{(1)}(n)$.}
\end{center}
\end{table}
Note that the expression for the first derivative obtained by this
method coincides with the value presented in the
Refs.~\cite{AbrSti64,Bic41} for $n=1,2$ that corresponds to three
and five points of interpolation. However, technique developed in
this article allows one to calculate the coefficients
$\alpha_{m}^{(1)}(n)$, and hence the first derivative, for any
value of $n$.

Similar formula can be obtained for the calculation of the second
derivative. We give without proof corresponding expression
\begin{equation}
  \label{f''}
  f^{\prime\prime}(0)\approx P_{2n}^{\prime\prime}(0)=
  \frac{1}{h^{2}}
  \sum_{m=1}^{n}\alpha_{m}^{(2)}(n)(f_{m}-2f(0)+f_{-m}),
\end{equation}
where
\begin{equation}
  \label{alpha2}
  \alpha_{m}^{(2)}(n)=\frac{1}{m^{2}\pi_{m}(n)},
\end{equation}
and the product $\pi_{m}(n)$ is introduced in Eq.~\eqref{pim}.

As an example of the use of the obtained central differentiation
formulae we will compute the first derivative of the function
$y(x)=\sin x$ at $x=0$. Let us set the value of the mesh width $h$
equal to $\pi/2$. Notice that, as a rule, the less the mesh width
$h$ the more exact result numerical differentiation gives. We have
chosen rather big value of $h$. The Eq.~\eqref{f'} for this case
takes the form
\begin{equation}
  \label{y'} y^{\prime}(0)= \frac{2}{\pi}
  \sum_{m=0}^{n}(-1)^{m}\alpha_{2m+1}^{(1)}(n).
\end{equation}
In Eq.~\eqref{y'} we take that $d_{m}=0$ if $m$ is an even number,
and $d_{m}=\pm 2$ if $m$ is an odd number.

Let us study the obtained result in the limiting case
$n\to\infty$. First, it is necessary to calculate the value of the
product $\pi_{m}(n)$ within the limit $n\to\infty$
\begin{multline}
  \label{limpi}
  \pi_{m}=\lim_{n\to\infty}\pi_{m}(n)=\lim_{\varepsilon\to 0}
  \frac{1}{1-
  \left(
  \frac{m+\varepsilon}{m}
  \right)^{2}
  }
  \prod_{k=1}^{\infty}
  \left(
  1-\frac{(m+\varepsilon)^{2}}{k^{2}}
  \right)=
  \\
  \frac{(-1)^{m+1}}{2}\lim_{\varepsilon\to 0}
  \frac{\sin\pi\varepsilon}{\pi\varepsilon}=
  \frac{(-1)^{m+1}}{2}
\end{multline}
Here we used the known value of infinite product
\[
  \prod_{k=1}^{\infty}
  \left(
  1-\frac{x^{2}}{k^{2}}
  \right)
  =\frac{\sin\pi x}{\pi x}.
\]
Using Eqs.~\eqref{alpha1} and \eqref{limpi}, we find that the
expression for the coefficients $\alpha_{m}^{(1)}(n)$ within the
limit $n\to\infty$ is represented in the following way
\begin{equation}
  \label{limalpha}
  \alpha_{m}^{(1)}=\lim_{n\to\infty}\alpha_{m}^{(1)}(n)=
  (-1)^{m+1}\frac{2}{m}.
\end{equation}
Now it is easy to complete the studying of Eq.~\eqref{y'}.
Substituting the result from Eq.~\eqref{limalpha} to
Eq.~\eqref{y'} and using the known value of infinite series we get
that
\[
  y^{\prime}(0)= \frac{4}{\pi}
  \sum_{m=0}^{\infty}\frac{(-1)^{m}}{2m+1}=1.
\]
Thus, it is shown that the method of the derivatives finding,
developed in this paper, gives for the function $y(x)=\sin x$ the
value of the first derivative which coincides with the exact
analytical one even at rather crude mesh width.

\section{Spectral characteristics of weight coefficients sequences}

\begin{figure}
  \centering
  \includegraphics[scale=.7]{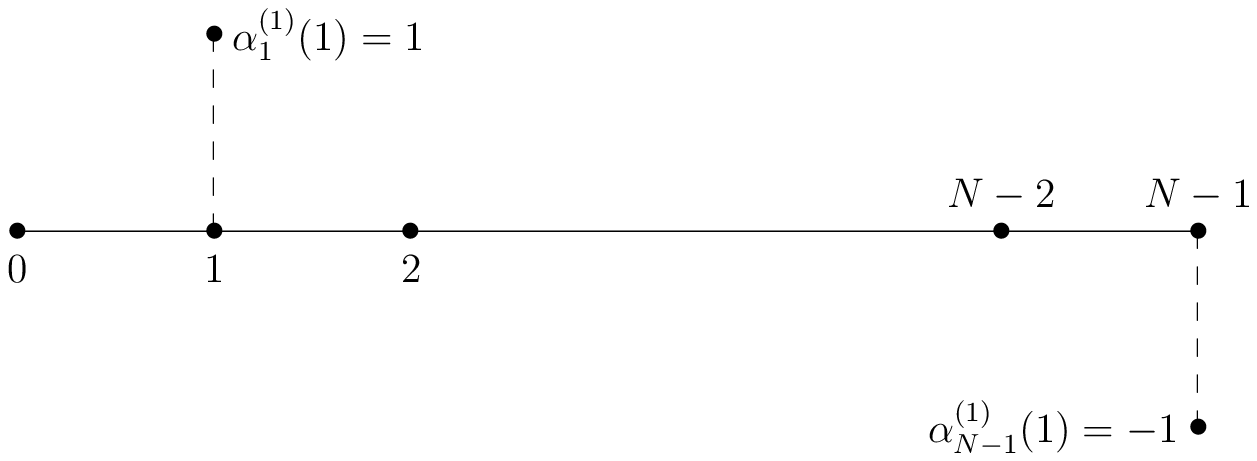}
  \caption{The coefficients $\alpha_{m}^{(1)}(1)$.}\label{alphan1}
\end{figure}

In the section~\ref{f'f''deriv} of the present work we derived the
formula for the finding of the first derivative of the function
$f(x)$ at $x=0$. This result can be easily generalized for the
case of the arbitrary point $x=kh$. If we set that
$\alpha_{-m}^{(1)}(n)=-\alpha_{m}^{(1)}(n)$, and moreover
supposing that $\alpha_{m}^{(1)}(n)=0$ for $m\equiv 0$ and $m>n$
(see Tab.~\ref{alphatab}), then in the considered case
Eq.~\eqref{f'} reads as follows
\begin{equation}
  \label{f'kh}
  f^{\prime}(kh)\approx P_{2n}^{\prime}(kh)=
  \frac{1}{2h}
  \sum_{m}\alpha_{m-k}^{(1)}(n)f_{m}.
\end{equation}
Here the summing is taken over all range of the function involved:
$f(kh)$. For instance, if the values of the function are set on
the limited equidistant collection of elements $N$, then
Eq.~\eqref{f'kh} can be rewritten in the form
\begin{equation}
  \label{f'khN}
  f^{\prime}(kh)\approx
  \frac{1}{2h}
  \sum_{m=0}^{N-1}\alpha_{m-k}^{(1)}(n)f_{m},
  \quad
  N\geq 2n.
\end{equation}
It is worth noticing that in Eq.~\eqref{f'khN} we used the
periodicity condition of the weight coefficients
\[
  \alpha_{m}^{(1)}(n)=\alpha_{m-N}^{(1)}(n).
\]
Fig.~\ref{alphan1} presents the example of the weight coefficients
of the differentiating sequence $\alpha_{m}^{(1)}(1)$.
Thus the first derivative computation of the function $f(x)$ at
the points $x=kh$, $k=0,1,\dots,N-1$, is reduced to the procedure
of the calculation of the mutual correlation function between the
finite sequences $\alpha_{m}^{(1)}(n)$ and $f_m$.

It is known (see, e.g., Ref.~\cite{Smi58p404eng}) that if a
function satisfies the Dirichlet conditions in the interval
$(-l,l)$, then it can be expanded into the Fourier series
\begin{equation}
  \label{fourierf}
  f(x)=\sum_{k=-\infty}^{+\infty}
  c_k\exp
  \left(
  i\frac{\pi k}{l}x
  \right),
\end{equation}
where the expansion coefficients are presented in the way
\[
  c_k=c^*_{-k}=
  \frac{1}{2l}
  \int_{-l}^{l}f(\xi)
  \exp
  \left(
  -i\frac{k\pi}{l}\xi
  \right)d\xi.
\]
If the first derivative $f'(x)$ satisfies the analogous conditions
as the function $f(x)$, then the following expression will be
valid
\begin{equation}
  \label{fourierf'}
  f'(x)=
  \sum_{k=-\infty}^{+\infty}
  \left\{
  i\frac{k\pi}{l}
  \right\}
  c_k
  \exp
  \left(
  i\frac{\pi k}{l}x
  \right),
\end{equation}
Therefore, form Eqs.~\eqref{fourierf} and \eqref{fourierf'} it
follows that the differentiation procedure is the linear filter
with the frequency characteristic: $\mathfrak{K}_1(k)=ik(\pi/l)$
\cite{KayMar81eng}. Similarly we receive for the second derivative
\[
  f''(x)=\sum_{k=-\infty}^{+\infty}
  \left\{
  -\left(
  \frac{k\pi}{l}
  \right)^2
  \right\}
  c_k
  \exp
  \left(
  i\frac{\pi k}{l}x
  \right),
\]
In this case the the frequency characteristic of the corresponding
filter has the form: $\mathfrak{K}_2(k)=-k^2(\pi/l)^2$.

According to Wiener-Khinchin theorem (see, e.g.,
Ref.~\cite{KayMar81eng}) the mutual correlation function between
the two finite sequences can be calculated with the help of the
inverse Fourier transform of the mutual spectrum of the considered
sequences. Thus, if we define that
\[
  \beta_1(r)=\sum_{m=0}^{N-1}
  \alpha_{m}^{(1)}(n)\exp
  \left(
  -i\frac{2\pi}{N}mr
  \right),
\]
is the complex spectrum of the differentiating sequence
$\alpha_{m}^{(1)}(n)$, and
\[
  c_r=\sum_{m=0}^{N-1}
  f_{m}\exp
  \left(
  -i\frac{2\pi}{N}mr
  \right),
\]
is the spectrum of the function $f(x)$, then it follows from
Eq.~\eqref{f'khN} that
\begin{equation}\label{mutcorr}
  f^{\prime}(kh)=
  \frac{1}{2h}
  \sum_{r=0}^{N-1}
  c_r\beta^*_1(r)
  \exp
  \left(
  i\frac{2\pi}{N}kr
  \right),
\end{equation}
where $\beta^*_1(r)$ is the complex conjugated quantity with
respect to $\beta_1(r)$.

Comparing Eqs.~\eqref{fourierf'} and \eqref{mutcorr} we obtain
that the accuracy of the numerical differentiation performed with
the use of the various types of the sequences
$\alpha_{m}^{(1)}(n)$ is characterized by the closeness of
imaginary parts of their spectra to the linearly growing sequence
$y_1(r)=2\pi r/N$.

\begin{figure}
  \centering
  \includegraphics[scale=.4]{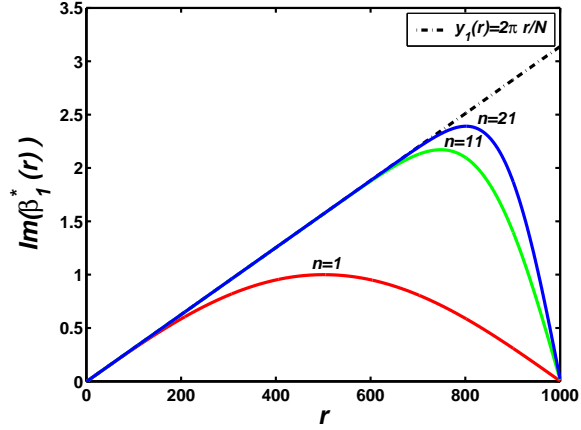}
  \caption{The spectra of various sequences
  $\alpha_{m}^{(1)}(n)$ at $N=2000$.
  }
  \label{spectra111211}
\end{figure}

The spectra of the sequences $\alpha_{m}^{(1)}(n)$ are depicted in
Fig.~\ref{spectra111211} for the various values of $n$ at
$N=2000$.
It can be seen from this figure that for $n=1$, i.e. for the
sequence shown in Fig.~\ref{alphan1}, the imaginary part of the
spectrum is the branch of the function $\sin(2\pi r/N)$. The
linearity condition is satisfied only in the vicinity of zero and
$N/2$. However, at $n=11$ the spectrum practically does not differ
from the linear one up to $r\approx N/2$. The more close to linear
one is the spectrum of the sequence $\alpha_{m}^{(1)}(21)$.

\begin{figure}
  \centering
  \includegraphics[scale=.4]{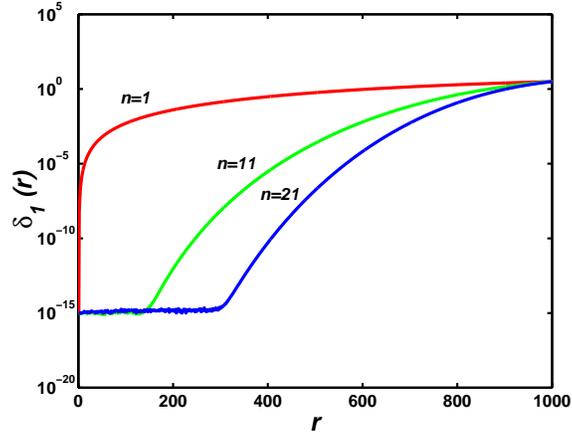}
  \caption{The function
  $\delta_1(r)=\Im\mathfrak{m}\big(\beta^*_1(r)\big)-y_1(r)$
  versus $r$ for different $n$.}
  \label{log101}
\end{figure}

The difference between the imaginary parts of the spectra of the
sequences $\alpha_{m}^{(1)}(n)$ and the linearly growing sequence
$y_1(r)=2\pi r/N$ are presented in Fig.~\ref{log101}.
The computations have been performed with the accuracy up to
$10^{-15}$, thus the reliable results at $n=11$ have been obtained
for $r\gtrsim 150$, and at $n=21$ for $r\gtrsim 300$. The
presented results demonstrate the high accuracy of the numerical
differentiation carried out with the help of the sequences
$\alpha_{m}^{(1)}(n)$ in the wide range of the spatial
frequencies.

\begin{figure}
  \centering
  \includegraphics[scale=.4]{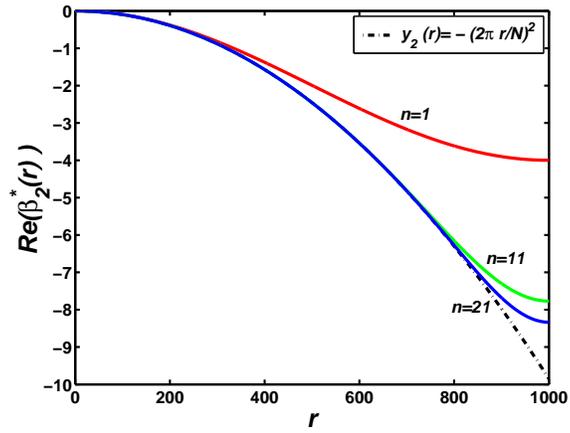}
  \caption{The spectra of various sequences
  $\alpha_{m}^{(2)}(n)$ at $N=2000$.
  }
  \label{spectra111212}
\end{figure}

\begin{figure}
  \centering
  \includegraphics[scale=.4]{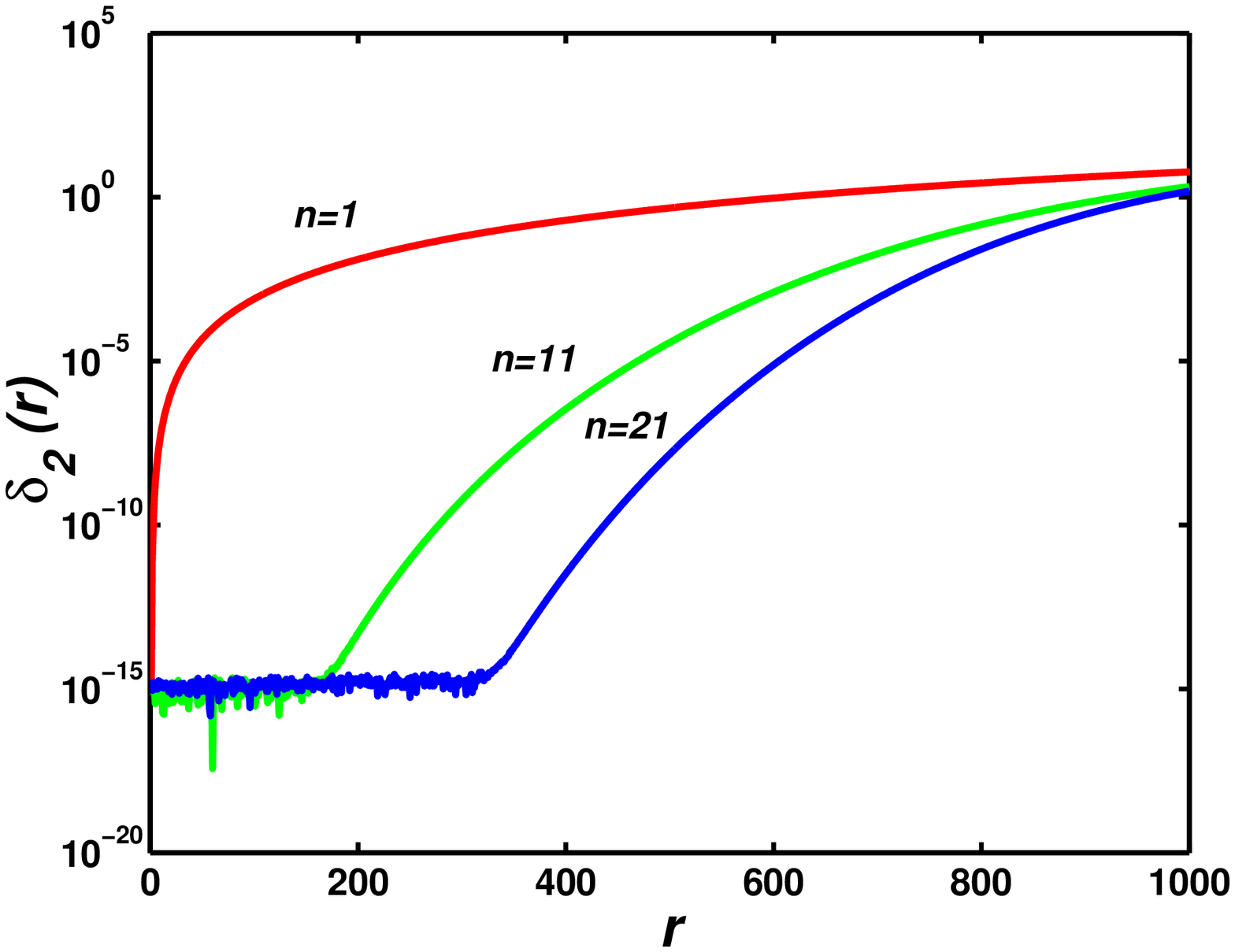}
  \caption{The function
  $\delta_2(r)=\Re\mathfrak{e}\big(\beta^*_2(r)\big)-y_2(r)$
  versus $r$ for different $n$.}
  \label{log102}
\end{figure}

Now let us briefly consider the sequences for the calculation of
the second derivative $\alpha_{m}^{(2)}(n)$, which are given in
Eq.~\eqref{alpha2}. Their spectral properties can be obtained in
the similar manner as we have done it for the case of the
sequences $\alpha_{m}^{(1)}(n)$ and therefore we just present the
final results. The spectra of the sequences $\alpha_{m}^{(2)}(n)$
are shown in Fig.~\ref{spectra111212}.
It follows from this figure that the closeness of the
corresponding spectrum to the parabola $y_2(r)=-(2\pi r/N)^2$ in
the case of $n=1$ exists only in the vicinity of zero. The spectra
at $n=11$ and $n=21$ are close to function $y_2(r)$ in a wider
range of $r$ ($r\lesssim 750$ and $r\lesssim 800$ respectively).
The difference between the real parts of the spectra of the
sequences $\alpha_{m}^{(2)}(n)$ and the parabola $y_2(r)=-(2\pi
r/N)^2$ are depicted in Fig.~\ref{log102} in the logarithmic
scale.
This figure again demonstrates the high accuracy of the second
derivative computation with the use of the sequences
$\alpha_{m}^{(2)}(n)$.

\section{Conclusion}

In conclusion we note that the method of central differentiation
formulae finding has been developed in this article. The
elaborated technique does not require direct construction of the
interpolating polynomial. We have derived simple and convenient
expressions for the first and the second derivatives
[Eqs.~\eqref{f'} and \eqref{f''}] of the function given in
$(2n+1)$ discrete points. The number $n$ was taken to be
arbitrary. In contrast to the results of the Ref.~\cite{For88},
where the recursion relations for the calculation of the weight
coefficient being used in numerical differentiation formulae were
considered, in the present work the expressions for the considered
weight coefficients have been derived in the explicit form for the
arbitrary number of interpolation points. As an example of the use
of the developed method we have calculated the first derivative of
the function $y(x)=\sin x$. The obtained result has been studied
in the limiting case $n\to\infty$. We have examined the spectral
characteristics of the weight coefficients sequence of the
numerical differentiation formulae for the different number of the
interpolation points. The performed analysis has allowed one to
study the accuracy of the numerical differentiation carried out
with the help of the developed method. It has been found that the
derived formulae of numerical differentiation posses the high
accuracy in a rather wide range of the spatial frequencies. As it
has been shown in this paper, the formulae for the derivatives
finding gave correct results in the case of large number of
interpolation points. Thus, the developed method can be useful in
lattice simulation of quantum fields \cite{Cre85}. To get the
exact results at calculations on lattices one has to use nets with
the big number of points. Derivatives which one encounters in
theories of quantum fields, as a rule, do not exceed the second
order. Therefore, the formulae obtained in this article could be
of use in carrying out mentioned above research.

\section*{Acknowledgments}

This research was supported by grant of Russian Science Support
Foundation. The author is indebted to Sergey~I.~Dvornikov for
helpful discussions.


\appendix

\section{Algorithm for computation of weight coefficients
$\alpha_{m}^{(1)}(n)$ \label{append}}

In this appendix we present the algorithm for the computation of
the coefficients $\alpha_{m}^{(1)}(n)$ on the MATLAB~6.5
programming language.

\begin{programbox}
 N=2000;\rcomment{$N$ is the size of the array}
 n=11;\rcomment{$n$ is the order of the differentiating sequence}
 \alpha=\mathtt{zeros}(1,N);\rcomment{$\alpha$ is the array of the weight coefficients}
 k_1=2; k_2=N;
    \FOR m=1:n
    r_1=1;
      \FOR k=1:n
      \IF k==m
      r_2=1;
      \ELSE
      r_2=1-(m/k)^2;
      \END
      r_1=r_1*r_2;
      \END
    r_1=1/(2*r_1*m);
    \alpha(1,k_1)=-r_1; k_1=k_1+1;
    \alpha(1,k_2)=r_1; k_2=k_2-1;
    \END
\end{programbox}

\end{document}